\documentclass{article}
\usepackage{amsmath,amsthm,amssymb}

\newcommand{\ga}{\alpha}
\newcommand{\gb}{\beta}
\renewcommand{\gg}{\gamma}
\newcommand{\gd}{\delta}
\newcommand{\gw}{\omega}

\newcommand{\gs}{\sigma}

\newcommand{\R}{\mathbb{R}}
\newcommand{\cP}{\mathcal{P}}

\newcommand{\dotrgen}{{\dot r}_{\mathit{gen}}}

\newcommand{\norm}{\mathit{norm}}
\newcommand{\Succ}{\mathit{succ}}

\newcommand{\cq}{\mathbb{Q}}

\newcommand{\Coll}[1]{\mathit{Coll}(\gw, <{#1})}

\newtheorem{theorem}{Theorem}[section]
\newtheorem{lemma}[theorem]{Lemma}
\newtheorem{claim}[theorem]{Claim}

\theoremstyle{definition}
\newtheorem{definition}[theorem]{Definition}
\newtheorem{example}[theorem]{Example}

\title{Games with creatures}
\author{{Saharon Shelah} and {Jind\v rich Zapletal}
\thanks{The first author is supported by The Israel Science 
Foundation funded by the Israel Academy of Sciences and Humanities. Publication 791.
The second author is partially supported by grants GA \v CR 201-00-1466, NSF DMS-0071437 
and a CLAS UF research award.}}

\begin{document}

\maketitle
\begin{abstract}
Many forcing notions obtained using the creature technology are naturally connected with certain integer games.
\end{abstract}

\section{Introduction}
The paper \cite{z:detandci} revealed an intimate connection between the proper forcing technology and determinacy
of infinite games. Many definable proper forcings $P$ adding a real turn out to have an integer game attached to them.
This means that there are Borel sets $C, D\subset\gw^\gw$ and a Borel function $f:C\to\R$ such that 

\begin{enumerate}
\item the set $I=\{A\subset\R:$ the second player
has a winning strategy in the integer game with payoff set $D\cup f^{-1}A\}$ is a $\gs$-ideal
\item every analytic $I$-positive set has a Borel $I$-positive subset
\item the poset $P$ is
forcing equivalent to the poset of all Borel $I$-positive sets ordered by inclusion.
\end{enumerate}

This is in fact a very natural concept, but it must be illustrated on examples.

\begin{example}
The game attached to Sacks forcing is just the perfect set game \cite{jech:set} Theorem 102. Just let $C=\gw^\gw$, $D=0$
and use a reasonable coding
of the perfect set game into an integer game to let $f:\gw^\gw\to\R$ be given by $f(x)=$the real
resulting from the run of the perfect set game coded by $x$.
\end{example}

\begin{example}
The game attached to Miller forcing is the unboundedness game of \cite{kechris:bounded}
\end{example}

\begin{example}
The game attached to Cohen forcing is the Baire category game \cite{jech:set} Theorem 102.
\end{example}

\begin{example}
Mathias forcing does not have an integer game attached to it. This follows from the fact that Mathias
forcing can be written as a two step $\gs$-closed*c.c.c. iteration, see the next paragraph.
\end{example}

In this paper we show that there is a pattern extending to many partial orders obtained by the creature technology of \cite{RoSh:470},
and in fact the proofs of properness of such posets naturally generate the related games. Why are we interested in the
connection? First of all, the concepts of both properness and determinacy have been around for several decades
and it is a priori interesting to see if they have anything in common. But there are practical reasons as well.
Under suitable determinacy assumptions every forcing $P$ with a game ideal $I$ attached to it satisfies a dichotomy:
every projective set of reals has either a Borel $I$-positive subset or a coanalytic $I$-small superset--Lemma~\ref{dich1}.
Such a dichotomy simplifies the treatment of the countable support iteration of the poset $P$ and the statement
of the absoluteness theorems in \cite{z:detandci}. Such a dichotomy also means that all
intermediate forcing extensions within the $P$-extension are c.c.c. \cite{z:closed}

Let us now describe the results of the paper.
The creature technology is quite complex and it is difficult--indeed self-defeating--to treat it in its full generality.
In the language of \cite{RoSh:470} we will deal with partial orders given by suitably definable tree-creating pairs of countable character
whose subcomposition operation is trivial, with no glueing. The information carried by such pairs can be
coded into a \emph{simple norm}, a Borel function whose domain is the set of all \emph{creatures} of the form
$\langle t, X\rangle$ where $t\in\gw^{<\gw}$ is a sequence and $X$ is a set of its one-step extensions. The values of
the function are nonnegative
real numbers. A simple norm gives rise to several forcing notions, among them

\begin{description}
\item[$\cq_0$] This is the partial order of all trees $T\subset\gw^{<\gw}$ such that for every natural number
$n$ and a node $s\in T$ there is a longer node $t\in T$ such that $\norm(t,\Succ_T(t))>n$. The ordering is by inclusion.
\item[$\cq_1$] This is the partial order of all trees $T\subset\gw^{<\gw}$ such that for every path $x\in[T]$
the limes inferior of the numbers $\norm(x\restriction n,\Succ_T(x\restriction n)):n\in\gw$ is infinity. The ordering is again by inclusion.
\end{description}

Of course these forcings depend on the choice of the norm, but in this paper all theorems are stated
for the general case of an arbitrary simple norm, and we will not mention the dependence. We
will always tacitly assume that the full tree $\gw^{<\gw}$ belongs to the poset $\cq_1$. 
There are many standard notions of forcing that can be obtained in this way: Sacks, Miller and Laver forcing, $PT_{fg}$ of
\cite{bartoszynski:set} 7.3.3, the $LT$ forcings of \cite{Sh:326} and others. 

\begin{example}
Let $\norm(t, X)$ be the length of the sequence $t$ if $X$ is infinite and $\norm(t, X)=0$ otherwise. Then $\cq_0$ is the Miller forcing
and $\cq_1$ is the Laver forcing.
\end{example}

The exact syntactical complexity of the norm is not critical as long as it is a projective function.  
For norms more complex than Borel we would have to change
the wording of our theorems a bit. However, the norms occurring in practice are very simple. Note that 
in the case of finitely branching trees
the norm is actually a hereditarily countable object.

It turns out that under our assumptions the forcings $\cq_0$ are always proper and they carry a particularly
simple game. This case is treated in Section 2. The forcings $\cq_1$ do not 
have to be proper, however the most common way of guaranteeing
the properness automatically provides an integer game for them. This is shown in Section 3.
In Section 4 we investigate the resulting dichotomies.

The notation of this paper follows the set theoretic standard of \cite{jech:set}. For a tree $T\subset\gw^{<\gw}$
and a node $t\in T$ the symbol $\Succ_T(t)$ denotes the set of immediate successors of the node $t$ in the tree $T$.
In forcing, we follow the western convention of writing $q\leq p$ when $q$ is a condition more informative than $p$. 
For forcings $\cq_0$ and $\cq_1$ the symbol $\dotrgen$ stands for the name for the
generic function in $\gw^\gw$, in both cases equal to the intersection of all
the trees in the generic filter. The abbreviation LC denotes a use of a suitable large cardinal hypothesis. 

\section{The $\cq_0$ forcing}
For the $\cq_0$ forcing arising from some simple norm, it is natural to consider the collection

\begin{definition}
$I_0$ is the collection of all sets $A\subset\gw^\gw$ which have an $F_\gs$-superset $B$ such that for no tree
$T\in\cq_1$ we have $[T]\subset B$.
\end{definition}

This collection can be restated in terms of an infinite game:

\begin{definition}
Let $A\subset\gw^\gw$ be a set. The infinite game $G_0(A)$ is played between Adam and Eve. First, Adam
plays a sequence $t_0\in\gw^{<\gw}$. After that, stage $n$ of the game proceeds as follows: a sequence
$t_n\in\gw^{<\gw}$ is given and Adam plays one by one some of its one step extensions (beginning a construction
of a creature of norm $\geq n$). The stage $n$ ends when (and if) Eve accepts one of these extensions, after
which Adam extends it further in an arbitrary way to a sequence $t_{n+1}\in\gw^{<\gw}$. Adam wins if either
at some stage $n$ Eve did not accept any of his extensions of the sequence $t_n$ and the creature he
constructed has norm $\geq n$, or else $\bigcup_nt_n\in A$.
\end{definition}

And the following is the key fact:

\begin{lemma}
\label{l1}
\begin{enumerate}
\item Adam has a winning strategy in the game $G_0(A)$ if and only if for some tree $T\in\cq_0$ it is the
case that $[T]\subset A$.
\item Eve has a winning strategy in the game
$G_0(A)$ if and only if $A\in I_0$.
\item $I_0$ is a $\gs$-ideal and $\cq_0$ is forcing equivalent to the algebra of Borel $I_0$-positive sets
ordered by inclusion.
\end{enumerate}
\end{lemma}

\begin{proof}
The first item is almost trivial. Suppose that Adam has a winning strategy in the game $G_0(A)$. Then the
initial segment closure of the set of all sequences that can possibly come up in a play following the strategy,
is some tree $T$, and it is not difficult to see that $T\in\cq_0$ and $[T]\subset A$. On the other hand, if
for some tree $T\in\cq_0$ we have $[T]\subset A$, then Adam can win the game $G(A)$ by playing only the nodes
that occur on the tree $T$. The reader can easily complete the standard argument.

For the second item, it is now clear from (1) and Borel determinacy that if $A\in I_0$
then Eve has a winning strategy in the game $G_0(A)$. For the converse, suppose that Eve has a winning strategy 
$\gs$ in the game $G_0(A)$ and fix some notation: given a partial run $v$ of the game which ends with Eve
accepting some sequence $t_v$ at a stage $n_v$, let $T_v=\{t\in\gw^{<\gw}:$ no intermediate node $t_v\subset
s\subseteq t$ can be accepted by Eve in some run of the game $G_0$ extending $v$ and respecting the strategy $\gs\}$.
It is not difficult to see that except for the initial segments of the sequence $t_v,$ the collection $T_v$
is closed under initial segment, and no node in it can branch into more than $\norm=n_v+1$ many immediate successors.
We will show that $A\subset\bigcup_v[T_v]$ and that no tree $T\in\cq_1$ has $[T]\subset\bigcup_v[T_v]$, 
concluding the proof of the second item. For the first statement, if $r\notin\bigcup_v[T_v]$ then it is possible
to construct runs $0=v_0\subset v_1\subset\dots$ respecting the strategy $\gs$ in such a way that
$t_{v_m}\subset r$ using the fact that at each stage $m$ of the construction, the real $r\in\gw^\gw$
must steer out of the tree $T_{v_m}$. Then $\bigcup_mv_m$ is a run of the game $G_0$ according to the
winning strategy $\gs$ resulting in the real $r$, meaning that $r\notin A$. For the second statement, choose a
tree $T\in\cq_0$. Enumerate all partial runs respecting
the strategy $\gs$ by $v_m:m\in\gw$ and construct nodes $0=t_0\subset t_1\subset$ so that the node $t_{m+1}$
branches into more than $\norm=n_{v_m}+1$ many immediate successors and the node $t_{m+2}$ falls out
of the tree $T_{v_m}$. Then $\bigcup_mt_m$ is a branch through the tree $T$ which is not in the set $\bigcup_v[T_v]$.

There are now two dual ways to show that the collection $I_0$ is a $\gs$-ideal; one uses a fusion argument with
the poset $\cq_0$ and the other will combine countably many Eve's strategies into one. Both of the arguments are
rather standard, and one point in this paper is that they are really the same thing. Let us perform the
fusion argument. Suppose that $A_n:n\in\gw$ are $F_\gs$ sets, each of them without a subset of the form
$[T]$ where $T\in\cq_0$. Clearly $\bigcup_nA_n$ is a $F_\gs$ set and it is enough to prove that it does not
contain all branches of some $\cq_0$ tree $T$. Suppose for contradiction it does; then $T\Vdash\dotrgen\in\bigcup_n\dot A_n$
and by thinning out the tree $T$ if necessary we may find a particular number $n$ such that $T\Vdash\dotrgen\in
\dot A_n$. Certainly $A_n$ is an analytic set and there is a tree $U\subset(\gw\times\gw)^{<\gw}$
which projects into $A_n$ in all generic extensions, in particular in the $\cq_0$ extension. So there
is a name $\tau$ for a function in $\gw^\gw$ such that $T\Vdash\langle\dotrgen,\tau\rangle$ is a branch through the
tree $\check U$. A standard fusion argument will now yield a condition $S\subset T$ in the poset $\cq_0$ such that
for every number $m$ the tree $S_m=\{s\in S: S\restriction s$ does not decide the value of $\tau(\check m)\}$
is wellfounded. This means that for every branch $x$ through the tree $S$ the expression $\tau/x$ makes sense,
it is a function in $\gw^\gw$ and $\langle x,\tau/x\rangle$ is a branch through the tree $U$. This is to say
that $[S]\subset A_n$, contradicting the properties of the set $A$.

It is now clear that the poset $\cq_0$ is isomorphic to the algebra of $I_0$-positive Borel sets ordered by
inclusion: the function $\pi(T)=[T]$ is an order isomorphism between the poset $\cq_0$ and a dense subset
of the algebra.
\end{proof}

Thus the poset $\cq_0$ has a game attached to it, namely the game $G_0$, and the corresponding ideal has
a basis consisting of $F_\gs$ sets. This has numerous consequences for the forcing properties of the poset, see
for example \cite{z:closed}.

Whenever $I$ is a $\gs$-ideal with a basis consisting of $F_\gs$ sets, then the algebra of $I$-positive Borel
sets is a proper notion of forcing \cite{z:closed}. There are very many such algebras and by far not all of
them are forcing equivalent to a poset of the form $\cq_0$, if only for the reason that closed sets are not
dense in them. A simple example is the ideal of meager sets: while the equivalence classes of clopen sets
are dense in the factor algebra of Borel sets modulo the meager sets, the closed sets themselves are not dense in the
original unfactored algebra of nonmeager Borel sets ordered by inclusion. A somewhat more sophisticated
example is provided by the ideal $\gs$-generated by closed measure zero sets. Clearly, no positive Borel measure zero
set can have a closed positive subset, and the situation carries over to the factor algebra as well.

\section{The $\cq_1$ forcing}
Now let us look at the forcing $\cq_1$ arising from some simple norm. It is natural to consider the following object:

\begin{definition}
$I_1$ is the collection of all sets $A\subset\gw^\gw$ which have a coanalytic superset $B\subset\gw^\gw$ such that
for no tree $T\in\cq_1$ it is the case that $[T]\subset B$.
\end{definition}

It is not difficult to express the collection $I_1$ in terms of an infinite game:

\begin{definition}
Let $A\subset\gw^\gw$ be a set. The infinite game $G_1(A)$ between players Adam and Eve
is played as follows. First Adam indicates a sequence $t_0\in\gw^{<\gw}$. Then, at stage $n$ of the game
a sequence $t_n$ will be known, Adam will play a natural number $m_n$ and then one-by-one a set of one-step extensions
of the sequence $t_n$ (beginning a construction of a creature of norm $>m_n$ based on $t_n$). 
The stage $n$ will be finished when (and if) Eve accepts one of these extensions
to become $t_{n+1}$. Adam wins if either at some stage $n$ Eve has not played for infinitely many rounds and
the creature Adam constructed at that stage has norm $>m_n$, or else the numbers $\{m_n:n\in\gw\}$ diverge to infinity
and the real $\bigcup_nt_n$ belongs to the set $A$. 
\end{definition}

We will use the occassion to fix some notation relevant to the game $G_1.$ The real
$\bigcup_nt_n$ will be referred to as the \emph{outcome} of the play. The numbers $m_n$ that Adam plays will be
called \emph{norms}. If $u\subset v$ are two runs of the game $G_1$ then we will say that the run $v$ 
$n$\emph{-extends} the run $u$
if all the norms played on $v\setminus u$ are greater or equal to $n$. A \emph{good} strategy for Eve is one that cannot
be defeated by reason of not accepting any node at some stage while Adam constructs a creature of a suitable
norm at that stage. Clearly, the set of all good strategies for Eve is a coanalytic set. A play of the game
is \emph{correct} if Eve accepts a node at all the infinitely many stages, and the norms played diverge to infinity. 

It is not difficult to verify

\begin{lemma}
\begin{enumerate}
\item Adam has a winning strategy in the game $G_1(A)$ if and only if there is a tree
$T\in\cq_1$ such that $[T]\subset A$.
\item (LC) Eve has a winning strategy in the game $G_1(A)\}$ if and only if $A\in I_1$.
\end{enumerate}
\end{lemma}

\begin{proof}
The first item is next to trivial. If, on one hand, there is a tree $T\in\cq_1$ such that $[T]\subset A$ then Adam can win the game $G_1(A)$
by simply making sure that $t_n\in T$, playing $m_n=\norm(t_n, \Succ_T(t_n))$, and on the $n$-th stage
simply enumerating all the immediate successors of the node $t_n$ in the tree $T$. On the other hand,
Adam has a winning strategy $\gs$ in the game $G_1(A)$ then the tree $T$ of all nodes that can occur
as some $t_n$ in some run of the game observing the strategy $\gs$, is an element of the poset
$\cq_1$ such that $[T]\subset A$, as the reader can easily verify.

For the second item, first prove the right-to-left inclusion. Suppose $\gs$ is Eve's winning strategy in the game $G_1(A)$. Then the set $A$ is included in the
coanalytic set $B=\gw^\gw\setminus$the set of all outcomes of the correct plays which observe the strategy $\gs$.
The strategy $\gs$ is still winning for Eve in the game $G_1(B)$, and by the first item the set $B$ does not
have a subset of the form $[T]$ for a tree $T\in\cq_1$. Thus $A\in I_1$. For the other inclusion, suppose
that $A\in I_1$ is a set. This means that $A$ is a subset of a coanalytic set $B$ without a $\cq_1$-tree in it. 
By the analytic determinacy and the first item, Eve has a winning strategy in the game $G_1(B)$, which is of course
winning even in the game $G_1(A)$.
\end{proof}

Comparing the situation with the parallel development in the previous section, now it would be natural to prove
that $I_1$ is a $\gs$-ideal; then, as in Lemma~\ref{l1}(3), the poset $\cq_1$ will be forcing equivalent to the algebra
of $I_1$-positive Borel sets ordered by inclusion and $G_1$ will be the game attached to the forcing
$P_1$. However, here the situation is complicated by the fact that the collection $I_1$ is not necessarily a $\gs$-ideal,
a problem that is connected with the possibility that the poset $\cq_1$ is not proper.

We will prove two implications: ``$\cq_1$ has continuous reading of names'' implies that ``$I_1$ is a $\gs$-ideal'',
which in turn implies that ``$\cq_1$ is proper''. Here, the continuous reading of names is the most common
tool to ensure the properness of the forcing $\cq_1$. Thus we can be satisfied to conclude that for the forcings
of the form $\cq_1$ the existence of an attached integer game is tightly connected with properness.
Note that if $I_1$ is a $\gs$-ideal then the generic filter on $\cq_1$ can be reconstructed from the
generic real as the set of all conditions in $\cq_1$ containing the real by Lemma 2.1 of \cite{z:detandci}.
This property can fail for certain variations of the tree creature forcings. 

It should be remarked that it is in general impossible to find a basis for the ideal $I_1$ consisting of
$F_\gs$ sets. The reason is that for a suitable norm the forcing $\cq_1$ can add a dominating
real, and such a real cannot be added by an algebra of Borel $I$-positive sets ordered by inclusion, where
$I$ is an ideal with basis consisting of $F_\gs$ sets. However, in many cases it is possible to find
a basis for the ideal $I_1$ consisting of rather simple Borel sets, for example the ideal
associated with the Laver forcing has a basis consisting of $G_\gd$ sets.

\begin{definition}
The forcing $\cq_1$ has \emph{continuous reading of names} if for every collection
$\{O_n:n\in\gw\}$ of open dense subsets of it and every tree $T\in\cq_1$ there is a subtree $S\subset T$ in
the poset $\cq_1$ such that for every number $n$ the set $S_n=\{t\in S:S\restriction t\notin O_n\}$ is wellfounded.
\end{definition}

While this notion has not been explicitely defined in print, it has very often been used.
Cf. Lemma 2.3.6(2) in \cite{RoSh:470}. It is not difficult to see 
that if the poset $\cq_1$ has continuous reading of names then it preserves $\aleph_1$ and
every real in the extension is an image of the generic real under a continuous function in the ground model--hence
the name.

\begin{lemma}
(LC) Suppose that the forcing $\cq_1$ has continuous reading of names. Then $I_1$ is a $\gs$-ideal.
\end{lemma}

\begin{proof}
Suppose $I_1$ is not a $\gs$-ideal; then there must be $I_1$-small sets $A_n:n\in\gw$
such that $\bigcup_nA_n\notin I$. Let $B_n$ be coanalytic sets such that $A_n\subset B_n$ and
for no tree $T\in\cq_1$ and no natural number $n$, $[T]\subset B_n$. Now
$\bigcup_nA_n\subset\bigcup_nB_n$ and $\bigcup_nB_n\notin I_1$ is a coanalytic set. By the analytic determinacy and
the previous lemma, there is a tree $T\in\cq_1$ such that $[T]\subset\bigcup_nB_n$. Clearly,
$T\Vdash\dotrgen\in\bigcup_n\dot B_n$ and strengthening the tree $T$ if necessary we may assume that there
is a specific number $n$ such that $T\Vdash\dotrgen\in\dot B_n$. Let $U\subset(\gw\times\gw_1)^{<\gw}$ be a tree
whose projection in all $\gw_1$-preserving extensions is the set $B_n$. Thus $T\Vdash\dotrgen\in p[\check U]$
and there must be a $\cq_1$-name $\tau$ for a function from $\gw$ to $\gw_1$ such that
$T\Vdash\langle\dotrgen,\tau\rangle$ is a branch through the tree $\check U$. By the continuous reading
of names, there is a subtree $S\subset T$ in the poset $\cq_1$ such that for every natural number $m$
the subtree $S_m=\{t\in S:S\restriction t$ does not decide the value of $\tau(\check m)\}$ is wellfounded.
This means that for every path $x\in p[S]$ the formula $\tau/x=\{\langle m,\ga\rangle:$ for some initial segment
$t\subset x$ with $S\restriction t\Vdash\tau(\check m)=\check\ga\}$ defines a total function from $\gw$ to
$\gw_1$ and the pair $\langle x,\tau/x\rangle$ constitutes a path through the tree $U.$ This is to say
that $[S]\subset B_n$ contradicting the choice of the set $B_n$.
\end{proof}

Note that we have not used anything concerning the definability of the ordering $\cq_1$ in the above proof.

\begin{lemma}
(LC) If $I_1$ is a $\gs$-ideal then the forcing $\cq_1$ is proper.
\end{lemma}

\begin{proof}
Suppose that $I_1$ is a $\gs$-ideal. Then the poset $\cq_1$ is forcing equivalent to the
algebra of Borel $I_1$-positive sets ordered by inclusion, and it is enough to argue
for the properness of the algebra. We will prove the following claim of independent interest:

\begin{claim}
\label{masterclaim}
(LC) Let $P$ be a poset adding a real $\dot r$ which is forced to fall out of all ground model coded
$I_1$-small sets. Let $M$ be a countable elementary submodel of a large enough structure. For every
condition $p\in P\cap M$ the set $\{\dot r/g:g\subset P\cap M$ is an $M$-generic filter$\}$ is $I_1$-positive.
\end{claim}

The lemma follows immediately from the claim applied to the special case of the algebra of Borel $I_1$
positive sets via Lemma 2.2. of
\cite{z:detandci}. To prove the claim, fix the poset $P$ and the name $\dot r$. Choose a measurable
cardinal $\kappa>|P|$ and fix the usual tree $U\subset(\gw\times\kappa)^{<\gw}$ projecting into the
coanalytic set of all good Eve's strategies for the game $G_1$.  Consider another infinite game $H$
with the following rules: player I first indicates a condition $p_0\in P$ and then produces one-by-one
open dense sets $O_n:n\in\gw$ of the poset $P$ and nodes $u_0\subset u_1\subset\dots$ in the tree $U$.
Thus $\bigcup_nu_n$ is a branch through the tree $U$ and its first coordinates constitute a good Eve's strategy
that we will denote by $\gs$. Meanwhile, player II is allowed to tread water--to wait for an arbitrary finite number
of rounds before his next move. He produces one-by-one conditions $p_0\geq p_1\geq\dots$ and partial runs
$0=v_0\subset v_1\subset v_2\subset\dots$ of the game $G_1$. Player II wins if the filter $g$ generated
by the conditions $\{ p_n:n\in\gw\}$ meets all the sets $\{ O_n:n\in\gw\}$, the run $\bigcup_n v_n$ is correct, it
respects the strategy $\gs$, and it results in the real $\dot r/g$.

\begin{claim}
\label{cc}
(LC) Player II has a winning strategy in the game $H$.
\end{claim}

Claim~\ref{masterclaim} follows. Fix Player II's winning strategy  $\tau$, let $M$ be a countable elementary
submodel of a large enough structure with $\tau\in M$, let $p\in P\cap M$ and consider the set $A=\{\dot r/g:p\in g, g\subset P\cap M$
is an $M$-generic filter$\}$. Suppose for contradiction that $A\in I_1$; then there must be Eve's good winning strategy
$\gs$ in the game $G_1(A)$. Let $N$ be an elementary submodel such that $M\cap\cP(P)=N\cap\cP(P)$ and the ordertype of the set $N\cap
\kappa$ is $\gw_1$; there is such a model due to the measurability of the cardinal $\kappa$. The tree $U\cap N$
still projects into the set of all Eve's good strategies, in particular, there is a branch through the
tree $U\cap N$ which gives the strategy $\gs$. Now simulate a run of the game $H$ in which player II follows his strategy
$\tau$ and player I puts $p=p_0$ and then enumerates the incriminated branch of the tree $U\cap N$ and all
the open dense subsets of the poset $P$ in the model $M$. Since Player II's moves come from the model $N$
and $\tau\in N$, player I's responses come from the model $N$ too (or $M$, which is the same thing due
to the choice of the model $N$). Let $g\subset P\cap M$ be the filter 
and let $v$ be the run of the game $G_1$ generated by player II's responses. Clearly, $\dot r/g\in A$
and $v$ is a run of the game $G_1(A)$ which respects the strategy $\gs$ which results in this real--that
is, Adam has won. Thus the strategy $\gs$ was not winning for Eve, contradiction.

The observation critical for the proof of Claim~\ref{cc} is that the payoff set of the game $H$ is Borel in the
space of all possible runs of the game, therefore the game is determined by \cite {martin:borel}. So it is enough
to derive a contradiction from the assumption that player I has a winning strategy $\tau$. Player II
will construct a counterplay consisting of conditions $p_n$ and partial runs $v_n$ using the
following induction hypothesis:

\begin{enumerate}
\item the condition $p_n$ and the run $v_n$ are played at the same time, $p_{n+1}\in O_n$
\item the partial run $v_n$ obeys the part of the strategy $\gs$ that is known by the time the run $v_n$
is played. The run $v_n$ ends with Eve accepting some sequence $t_n$ and $p_{n}\Vdash\check t_{n-1}\subset\dot r$. 
The run $v_{n+1}$ $n$-extends the run $v_{n}$.
\item The condition $p_n$ is $v_n, n$-good, meaning that for every good Eve's strategy $\gs$ which can output
the run $v_n$ there exists a condition $q\leq p_n$ such that $q\Vdash$ there is a complete run  of the game $G_1$
which $n$-extends the run $v_{n}$, it is correct, it respects the strategy $\gs$
and results in the real $\dot r$.
\end{enumerate}

The first two items say that if player II can perform the construction, the result will be a run of the game $H$
against the strategy $\tau$ in which he wins, contradicting the assumption on the strategy $\tau$. The third
item is just an extra induction hypothesis designed so as to make the inductive step possible.

In the beginning, the condition $p_0$ is given, and let $v_0=0$. The induction hypothesis is satisfied. Note that the
third item in this case is implied by the assumption that the real $\dot r$ is forced to fall out of all ground model
coded $I_1$-small sets. Suppose $p_n, v_n$ are known. To construct the condition $p_{n+1}$ and the run
$v_{n+1}$, first look at what happens if player II treads water from this point on. Of course, in this way he loses,
but the main point is that the strategy $\tau$ will have to continue playing, creating Eve's strategy $\gs$
and also, somewhere on the way, the open dense set $O_n\subset P$. Since the condition $p_n\in P$
is $v_n, n$-good, there is a condition $q\leq p_n$ with the properties described in the third item above. 
Strengthening the condition
$q$ if necessary, we may assume that $q\in O_n$. 

\begin{claim}
(LC) There is a partial
run $w$ of the game $G_1$ which $n$-extends the run $v_n$, in which Eve uses the strategy $\gs$, at the last move
in it she accepts some sequence $t_{n+1}$ properly extending $t_n$, and such that the condition $q$ is $w, n+1$-good.
\end{claim}

This will conclude the inductive step,
since player II will be able to play $p_{n+1}=q, v_{n+1}=w$ as soon as the strategy $\tau$ reveals
the open dense set $O_n$ and enough of the strategy $\gs$ to see that $w$ is a run in which Eve
uses that strategy, and the inductive hypothesis will continue to hold.

Well, suppose that the claim fails. Then for each run $w$ satisfying the
assumptions of the claim we can find a strategy $\gs_w$ for Eve in the game $G_1$
showing that the condition $q$ is not $w, n+1$-good. This means that $q\Vdash\dot r\in\dot X$ where
$X=\{s\in\gw^\gw:$ there is a complete correct run $n$-extending the run $v_n$, respecting the strategy $\gs$, and resulting
in the real $s$; however, for no partial run $w$ satisfying the assumptions of the claim there is a  correct complete
run $n+1$-extending the run $w$ which follows the strategy $\gs_w$ and results in the real $s\}$.
The definition of the set $X$ is a mouthful, but it is really a boolean combination of analytic sets.
By our large cardinal assumptions, the game $G_1(X)$ is determined and either $X\in I_1$ or $[T]\subset X$
for some tree $T\in\cq_1$. The first case is impossible, since $q\Vdash\dot r\in\dot X$ and the real
$\dot r$ is forced to fall out of all ground model coded $I_1$-small sets. So we are stuck in the second case.
There must be a node $t_{n+1}$ extending $t_n$ in the tree $T$ such that all the norms above that node are $\geq n+1$. There must be
a partial run $w$ $n$-extending the run $v_{n+1}$, in which Eve follows the strategy $\gs$ and in her last move
accepts exactly the sequence $t_{n+1}$. Now clearly Adam can $n+1$-extend the run $w$ into a correct complete run against
the strategy $\gs_w$ by playing only nodes of the tree $T$ below $t$. The result of that run should
be a real which is not in the set $X$ by its definition. At the same time, it will be a branch through the tree
$T$ and as such it does belong to the set $X$. Contradiction.
\end{proof}

\section{The dichotomies}

The natural conclusion of the previous sections are the following dichotomies. Consult \cite{z:detandci} for the 
definition of $I$-perfect sets and iterated Fubini power $I^\ga$ for $\gs$-ideals $I$ on the reals.

\begin{lemma}
\label{dich1}
(LC) Suppose that $P$ is a proper forcing with an integer game attached to it, as
witnessed by Borel sets $C,D\subset\gw^\gw$, a Borel function $f:C\to\R$ and the $\gs$-ideal $I$. 
For every universally Baire set $A\subset\gw^\gw$,

\begin{enumerate}
\item either $A$ has a Borel $I$-positive subset
\item or $A$ has a coanalytic $I$-small superset.
\end{enumerate}

For every countable ordinal $\ga$ and every universally Baire set $A\subset(\gw^\gw)^\ga$,

\begin{enumerate}
\item either $A$ has a Borel $I$-perfect subset
\item or $A$ is $I^\ga$-small.
\end{enumerate}
\end{lemma}

\begin{proof}
Suppose that $A\subset\gw^\gw$ is universally Baire set. Suitable large cardinal assumptions imply that the
integer game with the payoff set $D\cup f^{-1}A$ is determined. If player I has a winning strategy $\gs$, then
the set $f''\gs''\gw^\gw\subset A$ is analytic, it is $I$-positive since $\gs$ remains a winning strategy
for player I in the associated game, therefore it has a Borel $I$-positive subset and we are in the first
case of the dichotomy. If on the other hand player II has a winning strategy $\gs$ then the set
$\R\setminus f''\gs''\gw^\gw\supset A$ is coanalytic and it is in the ideal $I$ since the strategy
$\gs$ remains a winning strategy for player II in the associated game.

The second dichotomy is proved using the first dichotomy and the results of \cite{z:detandci}.
\end{proof}

\begin{lemma}
(LC) The above dichotomies hold for the forcing $\cq_0$. Suppose that the forcing $\cq_1$ has continuous reading of names. 
Then the above dichotomies hold for $\cq_1$ as well.
\end{lemma} 

All similar dichotomies beg a question:
what happens in the classical choiceless Solovay model?

\begin{lemma}
(LC) Suppose that the forcing $\cq_1$ has continuous reading of names and $\kappa$ is an inaccessible cardinal.
The above dichotomies hold in the model $V(\R)\subset V[G]$ for all subsets of $\gw^\gw$
and $(\gw^\gw)^\ga$ respectively, whenever $G\subset\Coll{\kappa}$ is a generic filter.
The same holds about the forcing $\cq_0$.
\end{lemma}

\begin{proof}
We will treat the case of the $\cq_1$ forcing, the $\cq_0$ case being much simpler.

Start with the first dichotomy. Suppose that $A\in V(\R)$ is a subset of $\gw^\gw$. 
By a standard argument we may assume that the set $A$ is definable in $V(\R)$ as $A=\{r\in\gw^\gw:\phi(r,t)\}$
from
some parameter $t$ in the ground model.
Assume that the set $A$ is forced to be $I_1$-positive. Then there must be in $V$ a forcing $P$ of size
$<\kappa$ and a $P$-name $\dot r$ such that $P\Vdash$``the real $\dot r$ falls out of all the $I_1$-small
sets coded in the ground model and $\Coll{\kappa}\Vdash V(\R)\models\phi(\dot r, \check t)$''.
By essentially Claim~\ref{masterclaim} applied in the model $V[G]$ to the model $V$ instead of the arbitrary countable $M$,
the set $B=\{r/g:g\subset P$ is $V$-generic$\}$ is $I_1$-positive. It is also a Borel set, and by the
choice of the name $\dot r$ and the homogeneity properties of the forcing $\Coll{\kappa}$, it is
also the case that $B\subset A$.

For the second dichotomy assume that $\ga\in\kappa$ is an ordinal and $A\subset(\gw^\gw)^\ga$ is a set in $V(\R)$.
Working in the model $V(\R)$, we will show that if the set $A$ is not $I_1^\ga$-small, then it contains an
$I_1$-perfect Borel subset. By a standard argument we may assume that $\ga\in\gw_1^V$ and the set $A$ is definable
as $A=\{\vec r\in(\gw^\gw)^\ga:\phi(\vec r, t)\}$ from some parameter $t$ in the ground model.  Consider the
following strategy $\gs$ for Adam in the transfinite game of length $\ga$ defining the ideal $I_1^\ga$:
the strategy $\gs$ applied to a string $\vec s$ of Eve's answers gives the set $\bigcup\{X:X\in I_1$ coded in the model
$V[\vec s]$. If the set $A$ is not $I_1^\ga$-small, 
there must be a sequence $\vec r\in A$ which is a legal counterplay against this strategy.
This means that back in $V$ there is a poset $P$ of size $<\kappa$ and a $P$-name $\vec r$ for an
$\ga$-sequence of reals such that $P\Vdash$``for every $\gb\in\check\ga$, $\vec r(\gb)\notin
\bigcup\{X:X\in I_1$ coded in the model $V[\vec r\restriction\gb]\}$; moreover $\Coll{\kappa}\Vdash\phi(\vec r,\check t)$''.

Back to the model $V(\R)$. Call a sequence $\vec s\in(\gw^\gw)^{\leq\ga}$ $P$\emph{-generic} if
there is a $V$-generic filter $g\subset P$ such that $\vec s\subset\vec r/g$. For such a sequence
$\vec s$ and a condition $p\in P$ we will say that $p$ is $\vec s$-good if there is a $V$-generic filter
$g\subset P$ containing the condition $p$ such that $\vec s=\vec r/g$. The following claim is
reminiscent of the classical preservation theorems for countable support iterations. Note that $P$ is a countable
set in the model $V(\R)$.

\begin{claim}
For every ordinal $\gb\leq \ga$, for every ordinal $\gg\in\gb$, every Borel $I_1$-perfect set
$C\subset(\gw^\gw)^\gg$ consisting of $P$-generic sequences, and every Borel function $f:C\to P$
such that for every sequence $\vec s\in B$ the condition $f(\vec s)$ is $\vec s$-good, there is a Borel $I_1$-perfect set
$B\subset(\gw^\gw)^\gb$ consisting of $P$-generic sequences such that $C=B\restriction\gg$ and for every sequence
$\vec s\in B$ the condition $f(\vec s\restriction\gg)$ is $\vec s$-good. 
\end{claim}

Once the claim has been proved, we will apply it with $\gb=\ga$ and $\gg=0$ to get a Borel $I_1$-perfect set
$B\subset(\gw^\gw)^\ga$ consisting of $P$-generic sequences. By a standard argument using the homogeneity of
the poset $\Coll{\kappa}$ we can conclude that $B\subset A$ and the second dichotomy follows.

The claim is proved by induction on the ordinal $\gb$. First the successor step. Suppose the induction hypothesis
holds at $\gb$ and we want to verify it at $\gb+1$ for some ordinal $\gg\leq\gb$, a Borel $I_1$-perfect set
$C\subset(\gw^\gw)^\gg$ and a Borel function $f:C\to P$. Use the induction hypothesis to get a 
Borel $I_1$-perfect set $D\subset
(\gw^\gw)^\gb$ consisting of $P$-generic sequences such that $C=D\restriction\gg$ and for every sequence 
$\vec s\in D$ the condition $f(\vec s\restriction\gg)$ is $\vec s$-good. Now consider the set
$B=\{\vec s\in(\gw^\gw)^{\gb+1}:\vec s\restriction\gb\in D$ and $\vec s$ is a $P$-generic sequence
such that the condition $f(\vec s\restriction\gg)$ is $\vec s$-good$\}$. The set $B$ is as required;
the only thing to verify is that for every sequence $\vec s\in D$ the set $\{r\in\gw^\gw:
\vec s^\smallfrown r\in B\}$ is not in the ideal $I_1$. But this follows essentially from Claim~\ref{masterclaim}
applied to the forcing $P/\vec s$ below the condition $f(\vec s\restriction\gg)$,
the name $\vec r(\gb)$ and the model $V[\vec s]$ in place of the model
$M$. Note that the real
$\vec r(\gb)$ is forced to fall out of all $I_1$-small sets in the model $V[\vec r\restriction\gb]$.

For the limit step, suppose that $\gb$ is a limit of an increasing sequence of ordinals $\langle\gb_n:n\in\gw\rangle$,
$\gg\in\gb_0$ is an ordinal and $C\subset(\gw^\gw)^\gg$ and $f:C\to P$ are objects as in the assumption
of the claim. Let $\langle O_n:n\in\gw\rangle$ be an enumeration of open dense subsets of the poset $P$ in the ground model $V$,
and by use the inductive hypothesis on the ordinals $\gb_n$ repeatedly to construct a sequence
$\langle B_n, f_n:n\in\gw\rangle$ such that $B_n\subset(\gw^\gw)^{\gb_n}$ are Borel $I_1$-perfect sets
and $f_n:B_n\to P$ are Borel functions such that $B_{n-1}=B_n\restriction\gb_{n_1}$
and for every sequence $\vec s\in B_n$ the condition $f_n(\vec s)$ is $\vec s$-good, it belongs to the
open dense set $O_n$ and it is smaller than $f_{n-1}(\vec s\restriction\gb_{n-1})$. Here it is understood
that $C=B_{-1}$ and $f=f_{-1}$. The construction is very easy to perform: at each number $n\in\gw$ first
apply the induction hypothesis at $\gb_n$ to get a set $B_n$ as asserted in the Claim, and then for every
sequence $\vec s\in B_n$ let $f_n(\vec r)$ be some condition in the open dense set $O_n$ smaller
than $f_{n-1}(\vec s\restriction\gb_{n-1})$ which is $\vec s$-good, say the first condition with this property
in some fixed enumeration of the poset $P$. In the end, let $B=\{\vec s\in(\gw^\gw)^\gb:\forall n\in\gw\ 
\vec s\restriction\gb_n\in B_n\}$. It is clear that the set $B$ is Borel $I_1$-perfect. Moreover,
every sequence $\vec s\in B$ is $P$-generic and the condition $f(\vec s\restriction\gg)$ is $\vec s$-good,
as witnessed by the $V$-generic filter $g\subset P$ obtained from the descending sequence
$\langle f_n(\vec s\restriction\gb_n):n\in\gw\rangle$ of conditions in the poset $P$.

\end{proof}

\def\germ{\frak} \def\scr{\cal} \ifx\documentclass\undefinedcs
  \def\bf{\fam\bffam\tenbf}\def\rm{\fam0\tenrm}\fi 
  \def\defaultdefine#1#2{\expandafter\ifx\csname#1\endcsname\relax
  \expandafter\def\csname#1\endcsname{#2}\fi} \defaultdefine{Bbb}{\bf}
  \defaultdefine{frak}{\bf} \defaultdefine{mathfrak}{\frak}
  \defaultdefine{mathbb}{\bf} \defaultdefine{mathcal}{\cal}
  \defaultdefine{beth}{BETH}\defaultdefine{cal}{\bf} \def\bbfI{{\Bbb I}}
  \def\mbox{\hbox} \def\text{\hbox} \def\om{\omega} \def\Cal#1{{\bf #1}}
  \def\pcf{pcf} \defaultdefine{cf}{cf} \defaultdefine{reals}{{\Bbb R}}
  \defaultdefine{real}{{\Bbb R}} \def\restriction{{|}} \def\club{CLUB}
  \def\w{\omega} \def\exist{\exists} \def\se{{\germ se}} \def\bb{{\bf b}}
  \def\equivalence{\equiv} \let\lt< \let\gt>

\end{document}